\documentclass{amsart}
\usepackage{amssymb,latexsym}
\usepackage{graphicx}
\numberwithin{equation}{section}

\newcommand{\nn}{\nonumber}
\newcommand{\speq}{\!\!\!\!&=&\!\!\!\!}
\newcommand{\e}{\mathrm{e}}

\begin{document}

\title{Complex-type numbers and generalizations of the Euler identity}
\author{D. Babusci, G. Dattoli, E. Di Palma,  E. Sabia} 

\address{INFN - Laboratori Nazionali di Frascati, via E. Fermi 40, I-00044 Frascati, Italy.}
\email{danilo.babusci@lnf.infn.it}

\address{ENEA - Centro Ricerche Frascati, via E. Fermi 45, I-00044 Frascati, Italy.}
 \email{giuseppe.dattoli@enea.it}

\address{ENEA - Centro Ricerche Frascati, via E. Fermi 45, I-00044 Frascati, Italy.}
 \email{emanuele.dipalma@enea.it}
                 
\address{ENEA - Centro Ricerche Frascati, via E. Fermi 45, I-00044 Frascati, Italy.}
\email{elio.sabia@enea.it}

\begin{abstract}
We consider different generalizations of the Euler formula and discuss the properties of the associated 
trigonometric functions. The problem is analyzed from different points of view and it is shown that it can 
be formulated in a natural way in algebraic and geometric terms.
\end{abstract}

\maketitle

\section{Introduction}\label{sec:intro}
The circular trigonometric functions have long been recognized as the characteristics functions associated with the exponentiation of 
imaginary numbers and more in general of  matrices. Within such a framework the notion of trigonometry is not unique and we can 
consider different ``trigonometries" whose common thread can be traced back to some underlying algebraic structure. This point of 
view has inspired the researches of Ref. \cite{Herr1}, where a systematic group theoretic derivation of different formulations of 
trigonometries has been undertaken. The need for a wider conception of trigonometry has been pointed out in Refs. \cite{Yaglo,BirNo}, 
where the trigonometry in Minkowskian spacetime has been studied and exposed in depth. 

Circular functions and complex numbers are intimately connected and the most representative link is provided by the Euler identity
\begin{equation}
\label{eq:eucir}
\e^{\,\imath\,\theta} = \cos\theta + \imath\,\sin\theta\,.
\end{equation}
The properties of the circular functions can be derived, using a purely algebraic point of view, from eq. \eqref{eq:eucir} and from the 
cyclical properties of the \emph{circular} imaginary unit $\imath$, satisfying the identity $\imath^2 = - 1$. Furthermore, the differential 
equations satisfied by the ordinary trigonometric functions are just a consequence of the fact that the exponential function is an eigenfunction 
of the derivative operator, while the addition formulae are a consequence of the semi-group property of the exponential.

The introduction of a \emph{hyperbolic} ``imaginary" unit such that $\jmath^2 = 1$ \cite{FjeGal}, allows the definition of the hyperbolic functions, 
which satisfy the identity
\begin{equation}
\label{eq:euhyp}
\e^{\,\jmath\,\theta} = \cosh\theta + \jmath\,\sinh\theta\,.
\end{equation}

The theory of the trigonometric functions (in their wider form, circular and hyperbolic) can be developed independently of any geometrical 
interpretation and, within such a context, the distinction between the units $\imath$ and $\jmath$ does not appear particular meaningful: 
both can be viewed as solution of a second degree algebraic equation. We can therefore consider the following realization of a more general 
``imaginary" unit \cite{Yama}
\begin{equation}
\label{eq:genun}
h^2 = a + b\,h \qquad\qquad (a, b \in \mathbb{R})\,.
\end{equation}
By interpreting this equation as a second degree algebraic equation, its solutions 
\begin{equation}
h_\pm = \frac{b \pm \sqrt{\Delta}}{2} \qquad \qquad (\Delta = b^2 + 4\,a)
\end{equation}
can be used to introduce the following identities
\begin{equation}
\label{eq:eulgen}
\e^{\,h_\pm\,\theta} = C (\theta) + h_\pm \,S (\theta)\,, 
\end{equation}
that can be viewed as a generalization of the Euler formula. We have indicated with $C$ and $S$ two functions of the variable $\theta$  
that play the role of the trigonometric functions in eqs. \eqref{eq:eucir} and \eqref{eq:euhyp}. For this reason we refer to them as 
trigonometric-like functions (TLF). It is worth stressing that in writing eq. \eqref{eq:eulgen} we have assumed that $h_\pm$ are conjugated 
imaginary units, i.e., $h_+^c = h_-$, as in the case of circular and hyperbolic units.

By combining the previous equations, we can identify the functions $C$ and $S$ as
\begin{equation}
\left(
\begin{array}{c}
      C (\theta)    \\
      S (\theta)   
\end{array}
\right) = \frac{1}{\sqrt{\Delta}}\,\left(
\begin{array}{cc}
   - h_-   &  h_+  \\
      1     &  - 1 
\end{array}
\right)\,\left(
\begin{array}{c}
      \e^{\,h_+\,\theta}    \\
      \e^{\,h_-\,\theta}
\end{array}
\right)\,. 
\end{equation}
The differential equations satisfied by this family of functions is obtained by keeping the derivative with respect to $\theta$ of both sides of 
eq. \eqref{eq:eulgen} with respect to $\theta$ and, by using eq. \eqref{eq:genun}, we end up with (we use $h$ without any index because 
the identity depends only on the property \eqref{eq:genun} of the new imaginary unit)
\begin{equation}
C^\prime (\theta) + h\,S^\prime (\theta) = h\,C(\theta) + (a + b\,h)\,S (\theta)\,, \nn
\end{equation}
that, after equating the coefficients of terms in $h$, yields the following first-order differential equation
\begin{equation}
\label{eq:diffeq}
C^\prime (\theta) = a\,S (\theta)\,\qquad\qquad S^\prime (\theta) = C(\theta) + b\,S (\theta)\,.
\end{equation}
Moreover, it is easy to obtain the addition formulae \cite{FjeGal2}
\begin{equation}
\label{eq:addit}
\left(
\begin{array}{c}
      C_{12}    \\
      S_{12}   
\end{array}
\right) = \left(
\begin{array}{cc}
   C_1  &  a\,S_1  \\
   S_1  &  C_1 + b\,S_1 
\end{array}
\right)\,\left(
\begin{array}{c}
      C_2    \\
      S_2   
\end{array}
\right)\,,
\end{equation}
where we put
\begin{equation}
 f_k = f (\theta_k)\, \qquad\qquad f_{kl} = f (\theta_k + \theta_l) \qquad (f \equiv C, S;\, k,l = 1,2)\,. \nn
\end{equation}

It is worth noting that the functions $C$ and $S$ derived from eq. \eqref{eq:eulgen} can be expressed in terms of ordinary trigonometric functions 
(circular or hyperbolic), and therefore they cannot be therefore considered a new set of functions. However, if considered independently from the 
existence of their classical counterparts, they provide a tool having global properties (e.g., the addition formulae and the differential equations 
they satisfy) which make them very useful in applications, as, for example, the evolution of multilevel quantum systems (as discussed in the 
concluding section), or to treat problems of analytical mechanics \cite{Yama,Yama2}. 

The following general result holds: \emph{Given $\alpha, \beta \in \mathbb{R}$, it is always possible build a set of cosine- and sine-like functions 
according to the following prescription}
\begin{equation}
\label{eq:cossin}
C (\theta) = \frac{\alpha\,\e^{\,\beta\,\theta} - \beta\,\e^{\,\alpha\,\theta}}{\alpha - \beta}\,\qquad
S (\theta) = \frac{\e^{\,\alpha\,\theta} - \e^{\,\beta\,\theta}}{\alpha - \beta}\,,
\end{equation}  
whose proof is straightforward and is achieved by treating $\alpha, \beta$ as the independent roots of a second degree equation.

\section{Trigonometries and $2 \times 2$ matrices}\label{sec:trigo}
The notion of imaginary units exploited in the previous section can be made less abstract by noting that both circular and hyperbolic imaginary 
units can be expressed in terms of  matrices, namely\footnote{In the following, a superimposed ``hat"  denotes a matrix.}
\begin{equation}
\hat{\imath} = \left(
\begin{array}{cc}
   0   & -1   \\
    1  &  0 
\end{array}
\right)\,, \qquad\qquad
\hat{\jmath} = \left(
\begin{array}{cc}
   0   &  1   \\
    1  &  0 
\end{array}
\right)\,, 
\end{equation}
in such a way that  
\begin{equation}
\hat{\imath}^2 = -\,\hat{1}\,,\qquad\qquad \hat{\jmath}^2 = \hat{1}\,,
\end{equation}
with $\hat{1}$ the unit $2 \times 2$ matrix. As well known, the exponentiation of  $\hat{\imath}$ (the symplectic matrix) can be viewed as the 
generator of rotations in a plane, while the exponentiation of $\hat{\jmath}$ can be associated with non Euclidean rotations, and, more in 
general, with the Lorentz transformations viewed as rotations in the Minkowski plane.

By introducing the vector $\mathbf{z} (\theta) = (C (\theta), S (\theta))^T$, eq. \eqref{eq:diffeq} can be rewritten as follows
\begin{equation}
\label{eq:derCS}
\frac{\mathrm{d}}{\mathrm{d} \theta}\,\mathbf{z} (\theta) = \hat{h}\,\mathbf{z} (\theta)\,, \qquad\qquad 
\hat{h} = \left(
\begin{array}{cc}
   0 & a   \\
   1 & b 
\end{array}
\right)\;,
\end{equation}
with the matrix $\hat{h}$ satisfying the identity \eqref{eq:genun}. According to this equation the TLF are associated with the non-unitary 
``evolution operator"
\begin{equation}
\label{eq:umatr}
\hat{U} (\theta) = \e^{\,\theta\,\hat{h}}\,,
\end{equation}
The use of standard algebraic methods 
\cite{Datto} allows to derive the following expressions for the entries of the matrix associated to this operator
\begin{eqnarray}
U_{11} (\theta) \speq \left\{- \frac{b}{\sqrt{\Delta}}\,\sinh\left(\frac{\sqrt{\Delta}}{2}\,\theta\right) + 
\cosh\left(\frac{\sqrt{\Delta}}{2}\,\theta\right)\right\}\,\exp\left(\frac{b}2\,\theta\right) \nn \\
U_{22} (\theta) \speq \left\{\frac{b}{\sqrt{\Delta}}\,\sinh\left(\frac{\sqrt{\Delta}}{2}\,\theta\right) + 
\cosh\left(\frac{\sqrt{\Delta}}{2}\,\theta\right)\right\}\,\exp\left(\frac{b}2\,\theta\right)  \\ 
\frac{U_{12} (\theta)}{a} \speq U_{21} (\theta) = \frac2{\sqrt{\Delta}}\,\sinh\left(\frac{\sqrt{\Delta}}{2}\,\theta\right)\,
\exp\left(\frac{b}2\,\theta\right)\,. \nn
\end{eqnarray}
It is worth noting that the functions $C$ and $S$ can be viewed as independent linear combinations of the entries of the 
matrix \eqref{eq:umatr}, or, what is the same, that, according to eq. \eqref{eq:addit}, 
the evolution matrix can be written as
\begin{equation}
\hat{U} (\theta) = \left(
\begin{array}{cc}
   C (\theta)   &   a\,S (\theta) \\
   S (\theta)   &   C (\theta) + b\,S (\theta)
\end{array}
\right)\,.
\end{equation}
The previous results can be extended to any generic $2 \times 2$ matrix
\begin{equation}
\hat{m} = \left(
\begin{array}{cc}
   a   &  b  \\
   c   &  d 
\end{array}
\right)\,.
\end{equation}
In fact, the Cayley-Hamilton theorem \cite{FuOik} allows to prove the identity
\begin{equation}
\label{ }
\hat{m}^2 = - (\det \hat{m})\,\hat{1} + (\mathrm{tr}\,\hat{m})\,\hat{m}
\end{equation} 
that is a generalization of eq. \eqref{eq:genun}, and, therefore, all the functions that satisfy the differential equation
\begin{equation}
\label{eq:Udiff}
\frac{\mathrm{d}}{\mathrm{d} \theta}\,\mathbf{z} (\theta) = \hat{m}\,\mathbf{z} (\theta)\,,
\end{equation}
are TLF. The use of eq. \eqref{eq:genun} with $\hat{m}$ in place of $\hat{h}$, allows to write the entries of the evolution matrix 
associated to eq. \eqref{eq:Udiff} as follows
\begin{equation}
\label{eq:Uoper}
\hat{U} (\theta) = \e^{\,\theta\,\hat{m}} = \left(
\begin{array}{cc}
  C (\theta) + a\,S (\theta)  &  b\,S (\theta)  \\
  c\,S (\theta)  &  C (\theta) + d\,S (\theta)
\end{array}
\right)\,,
\end{equation}
that, compared with the expression of the matrix entries algebrically obtained 
\begin{eqnarray}
U_{11} (\theta) \speq \left\{\frac{a - d}{\sqrt{\Lambda}}\,\sinh\left(\frac{\sqrt{\Lambda}}{2}\,\theta\right) + 
\cosh\left(\frac{\sqrt{\Lambda}}{2}\,\theta\right)\right\}\,\exp\left(\frac{a + d}2\,\theta\right) \nn \\
U_{22} (\theta) \speq \left\{- \frac{a - d}{\sqrt{\Lambda}}\,\sinh\left(\frac{\sqrt{\Lambda}}{2}\,\theta\right) + 
\cosh\left(\frac{\sqrt{\Lambda}}{2}\,\theta\right)\right\}\,\exp\left(\frac{a + d}2\,\theta\right)  \nn \\ 
\frac{U_{12} (\theta)}{b} \speq \frac{U_{21}}{c} (\theta) = \frac2{\sqrt{\Lambda}}\,\sinh\left(\frac{\sqrt{\Lambda}}{2}\,\theta\right)\,
\exp\left(\frac{a + d}2\,\theta\right) \nn \\
& & \quad \left(\Lambda = (\mathrm{tr}\,\hat{m})^2 - 4\,\det \hat{m} = (a - d)^2 + 4\,b\,c\right)\,, \nn 
\end{eqnarray}
allows to obtain the explicit form for these TLF. As a consequence of eq. \eqref{eq:Uoper}, the following identity holds
\begin{equation}
\label{eq:CSide}
\exp\left(\theta\,\mathrm{tr}\,\hat{m}\right) = C^2 (\theta) + (\mathrm{tr}\,\hat{m})\,C (\theta)\,S (\theta) + 
(\det \hat{m})\,S^2 (\theta)\,,
\end{equation}
whose relevant geometrical and physical meaning will be discussed in the concluding section.

\section{Higher order trigonometries}\label{sec:hgtrig}
According to the analysis developed so far, what we have called trigonometry is the by-product of the formalism associated with 
the exponentiation of $2 \times 2$ matrices. ÒHigher order trigonometriesÓ can accordingly be considered along with the 
exponentiation of matrices with larger dimensionalities.

By following Ref. \cite{Yama2} we propose a further generalization of the Euler identity, by setting\footnote{We adopt the symbols 
$A_j (j = 0, 1, 2)$ for the trigonometric-like functions because more convenient and easily amenable for the $n$-dimensional 
extensions.}
\begin{equation}
\label{eq:exeta}
\e^{\,\eta\,\theta} = A_0 (\theta) + \eta\,A_1 (\theta) + \eta^2\,A_2 (\theta)\,,\qquad
\eta^3 = a_0 + \eta\,a_1 + \eta^2\,a_2\,.
\end{equation}
The expression of the functions $A_j$ can be found adopting exactly the same procedure as before, albeit slightly 
algebraically more cumbersome.
By setting
\begin{equation}
\mathbf{A} (\theta) = \left(A_0 (\theta),  A_1 (\theta), A_2 (\theta)\right)^T\,,\qquad\qquad 
\mathbf{Y} (\theta) = \left(\e^{\,\eta_1\,\theta}, \e^{\,\eta_2\,\theta}, \e^{\,\eta_3\,\theta}\right)\,, \nn
\end{equation} 
with $\eta_j\;(j = 1,2,3)$ roots of the cubic equation in eq. \eqref{eq:exeta}, we obtain
\begin{equation}
\label{eq:Afunct}
\mathbf{A} (\theta) = \hat{V}^{- 1}\,\mathbf{Y} (\theta)\,,
\end{equation}
where $\hat{V}$ is the Vandermonde matrix
\begin{equation}
\hat{V} = \left(
\begin{array}{ccc}
 1  &  \eta_1  &   \eta_1^2  \\
 1  &  \eta_2  &   \eta_2^2  \\
 1  &  \eta_3  &   \eta_3^2  
\end{array}
\right)\;. \nn
\end{equation}
Since a triple of numbers can always be associated with a cubic equation, it is clear that we can construct a set of third order 
TLF for each of these triple following the prescription \eqref{eq:Afunct}.

The differential equations satisfied by the third order TLF are derived by means of the same procedure as before, which yields
\begin{equation}
\frac{\mathrm{d}}{\mathrm{d} \theta} \mathbf{A} (\theta) = \hat{\eta}\,\mathbf{A} (\theta)
\end{equation}
with
\begin{equation}
\label{eq:lamatr}
\hat{\eta} = \left(
\begin{array}{ccc}
  0  &  0  &  a_0   \\
  1  &  0  &  a_1   \\
  0  &  1  &  a_2
\end{array}
\right)\,.
\end{equation}
The evolution matrix associated to the previous Cauchy problem can be written as 
\begin{equation}
\label{eq:eveta}
\hat{U} (\theta) = \e^{\,\theta\,\hat{\eta}} = A_0 (\theta) +  \hat{\eta}\,A_1 (\theta) + {\hat{\eta}}^2\,A_2 (\theta)\,,
\end{equation}
and, in terms of the functions $A_j$, reads
\begin{equation}
\hat{U} (\theta) = \left(
\begin{array}{ccc}
   A_0   &  a_0\,A_2  &  a_0\,A_1  + a_0\,a_2\,A_2   \\
   A_1   &  A_0  +  a_1\,A_2   &  a_1\,A_1  + (a_0 + a_1\,a_2)\,A_2  \\
   A_2   &  A_1  +  a_2\,A_2   &  A_0  +  a_2\,A_1  + (a_1 + a_2^2)\,A_2 
\end{array}
\right)\,.
\end{equation}
This equation can be used to compute the determinant of the evolution matrix in terms of the two-variable TLF, and from the relation
\begin{equation}
\det \hat{U} = \exp(\theta\,\mathrm{tr}\,\hat{\eta})
\end{equation}
we obtain
\begin{eqnarray}
\label{eq:3rdtrig}
\e^{\,a_2\,\theta} \speq A_0^3 + a_0\,A_1^3 + a_1\,A_2^3 - (3\,a_0 + a_1\,a_2)\,A_0A_1\,A_2 \nn \\
& & \;\;\;+ \,(2\,a_1 + a_2^2)\,A_2\,A_0^2 + (a_1^2 + 2\,a_0\,a_2)\,A_0\,A_2^2   \\
& &  \;\;\;- \,a_1\,A_1^2\,A_0 + a_0\,a_2\,A_1^2\,A_0 - a_0\,a_1\,A_1\,A_2^2 \nn
\end{eqnarray}
that should be viewed as the fundamental identity for the third order trigonometry.

As a particular case of third order TLF, let us consider the case of the Eisenstein numbers \cite{Gaal}.  They satisfy the properties 
\begin{equation}
\label{eq:eisen}
\omega^2 + \omega = - 1\,\qquad \omega^3 = 1\,,
\end{equation}
i.e., can be obtained letting $a_0 = 0$, $a_1 = a_2 = - 1$ and $\eta^3 = 1$ in the second of eq. \eqref{eq:exeta}. Following the treatment of 
sec. \eqref{sec:intro}, the associated $C$ and $S$ functions are easily obtained
\begin{eqnarray}
\label{eq:eitlf}
C (\theta) \speq \e^{\,-\,\theta/2}\,\left[\cos\left(\frac{\sqrt{3}}2\,\theta\right) + 
\frac1{\sqrt{3}}\,\sin\left(\frac{\sqrt{3}}2\,\theta\right)\right]\,, \nn \\
S (\theta) \speq \frac{2\,\e^{\,-\,\theta/2}}{\sqrt{3}}\,\sin\left(\frac{\sqrt{3}}2\,\theta\right)\,.
\end{eqnarray}
It is interesting to link the above functions to the so called pseudo-hyperbolic functions (PHF), discussed on the eve of the seventies of the last 
century \cite{Ricci}, and more recently discussed within the framework of an algebraic formalism in Refs. \cite{Datto1} and \cite{Datto2}. According 
to eqs. \eqref{eq:eisen} we find that
\begin{equation}
\label{eq:expeis}
\e^{\,\omega\,\theta} = e_0 (\theta) + \omega\,e_1 (\theta) + \omega^2\,e_2 (\theta)
\end{equation}
with
\begin{equation}
e_k (\theta) = \sum_{n = 0}^\infty \frac{\theta^{3\,n + k}}{(3\,n + k)!} \qquad\qquad (k = 0, 1, 2)\,. \nn
\end{equation}
Furthermore, from the condition
\begin{equation}
\e^{\,\omega\,\theta} = C (\theta) + \omega\,S (\theta)\,,
\end{equation}
by using eqs. \eqref{eq:eisen} and \eqref{eq:expeis}, we obtain
\begin{equation}
C (\theta) = e_0 (\theta) - e_2 (\theta)\,, \qquad\qquad S (\theta) = e_1 (\theta) - e_2 (\theta)\,.
\end{equation}

The extension of these results to the $n$-dimensional case is straightforward and computationally simplified by the remark that the 
structure of the matrices \eqref{eq:derCS} and \eqref{eq:lamatr} is that of a \emph{companion matrix} \cite{FuOik}. We remind that, 
given a $n \times n$ matrix $\hat{Q}$ with characteristic polynomials
\begin{equation}
P (\lambda) = \sum_{k = 0}^n p_{n - k}\,\lambda^k\,,\qquad\qquad (p_0 = 1)\,,
\end{equation}
its companion matrix is defined as
\begin{equation}
\label{eq:compa}
\hat{L} = \left(
\begin{array}{ccccl}
   0   & 0 & \cdots & 0 & - p_n \\
   1   & 0 & \cdots & 0 & - p_{n - 1} \\
   0   & 1 & \cdots & 0 & - p_{n - 2} \\
   \vdots & \vdots & \ddots & \vdots\\
   0   & 0 & \cdots & 1 & - p_1 
\end{array}
\right)
\end{equation}
and has the same eigenvalues of  $\hat{Q}$. Introduced the $n$-dimensional vector $\mathbf{e} = (1, 0, \cdots, 0)^T$, 
it is easy to show \cite{FuOik} that the following identity holds 
\begin{equation}
\hat{Q}^m = \left(\hat{1}, \hat{Q}, \cdots, \hat{Q}^{n - 1}\right)\,\left(\hat{L}^m\,\mathbf{e}\right) 
= \sum_{k = 0}^{n - 1} \hat{Q}^k\,\left(\hat{L}^m\,\mathbf{e}\right)_k\,.
\end{equation}
This result allows to conclude that
\begin{equation}
\label{eq:fuji}
\e^{\,\theta\,\hat{Q}} = \left(\hat{1}, \hat{Q}, \cdots, \hat{Q}^{n - 1}\right)\,
\left(\e^{\,\theta\,\hat{L}}\,\mathbf{e}\right)\,,
\end{equation}
and, since $\e^{\,\theta\,\hat{L}}$ can always be written in terms of $n$ suitable TLF, we can state that the exponentiation of any diagonalizable 
$n \times n$ matrix can be written in term of these functions.

Before concluding this section, let us come back to eq. (\ref{eq:eveta}) and consider its extension to the case of two variables 
\begin{equation}
\label{eq:ev2var}
\hat{U} (\theta, \phi) = \e^{\,\theta\,\hat{\eta} + \phi\,\hat{\eta}^2} = 
A_0 (\theta, \phi) +  \hat{\eta}\,A_1 (\theta, \phi) + {\hat{\eta}}^2\,A_2 (\theta, \phi)\,,
\end{equation}
Such an extension implies that the functions $A_k$ are actually two variable functions. Their forms are obtained quite straightforwardly 
by noting that the ``evolution" matrix (\ref{eq:ev2var}) and $\hat{\eta}$ are diagonalized by the same transformation. The functions 
$A_k (\theta, \phi)$ are therefore obtained by performing the substitution 
$\e^{\,\eta_j\,\theta}\,\to\,\e^{\,\eta_j\,\theta + \eta_j^2\,\phi}$ in eq. (\ref{eq:Afunct}). Let us note that in this case, 
the determinant of the matrix $\hat{U} (\theta, \phi)$ is given by
\begin{equation}
\det \hat{U} (\theta, \phi) = \exp\left(\theta\,\mathrm{tr}\,\hat{\eta} + \phi\,\mathrm{tr}\,\hat{\eta}^2\right)
\end{equation}
where 
\begin{equation}
\mathrm{tr}\,\hat{\eta} = \sum_{j = 1}^3 \eta_j = a_2\,,\qquad\qquad
\mathrm{tr}\,\hat{\eta}^2 = \sum_{j = 1}^3 \eta_j^2 = a_2^2 + 2\,a_1\,.
\end{equation}
and therefore the only modification in eq. (\ref{eq:3rdtrig}) consists in the multiplication of its lhs by the term 
$\e^{\,(a_2^2 + 2\,a_1)\,\phi}$.

\section{Concluding remarks}\label{sec:conclu}
As the usual trigonometric functions (circular or hyperbolic), also the TLFs have a geometric interpretation. Let us consider the identity given in 
eq. (\ref{eq:CSide}). In terms of the variables 
\begin{equation}
x  = \exp\left(-\,\frac{\theta}2\,\mathrm{tr}\,\hat{m}\right)\,C(\theta)\;,\qquad\qquad
y  = \exp\left(-\,\frac{\theta}2\,\mathrm{tr}\,\hat{m}\right)\,S(\theta)\;, \nn
\end{equation}
it can be rewritten as follows
\begin{equation}
\label{eq:quadra}
x^2  + (\mathrm{tr}\,\hat{m})\,x\,y + (\det \hat{m})\,y^2 = 1\,.
\end{equation}
The determinant of this quadratic form is
\begin{equation}
\delta = \frac{(\mathrm{tr}\,\hat{m})^2}4 - (\det \hat{m}) = \frac{(a - d)^2}4 + b\,c
\end{equation}
and its sign determine the geometric interpretation. For $\delta > 0$ the points with coordinates $x, y$, lie on an arc of hyperbola, and reduce to the 
ordinary hyperbolic functions when referred to the axes rotated with respect to $x, y$ by an angle
\begin{equation}
\chi = \frac12\,\tan^{- 1}\left(\frac{a + d}{a\,d - b\,c - 1}\right)\,.
\end{equation}
In the case $\delta < 0$ the quadratic form is represented by an ellipse and the same axis rotation reduce these functions to the 
circular ones. Finally, for $\delta = 0$ the $x, y$ functions does not make too much sense because the conic 
is degenerate.

As an example of third order TLF, we consider the particularly important case in which the matrix $\hat{\eta}$ has the following form 
\begin{equation}
\label{eq:neweta}
\hat{\eta} = \left(
\begin{array}{ccc}
  0  &  - \nu_3  &  \nu_2 \\
  \nu_3   &  0  &  - \nu_1 \\
  - \nu_2  &  \nu_1  & 0 
\end{array}
\right)\;.
\end{equation}
This matrix is not of the form (\ref{eq:lamatr}) and therefore the associated functions $A_k$ cannot be considered canonical. 
By using eq. (\ref{eq:Afunct}) it is easy to show that\footnote{$\mathrm{sinc}\,x = \sin x/x$.}
\begin{equation}
\label{eq:chaeta}
\e^{\,\theta\,\hat{\eta}} = \hat{1} + \theta\,\mathrm{sinc}\,(\nu\,\theta)\,\hat{\eta} + 
\frac12\,\theta^2\,\mathrm{sinc}^2\,(\frac{\nu}2\,\theta)\,\hat{\eta}^2 \qquad\
\left(\nu^2 = \sum_{i = 1}^3 \nu_i^2\right)\,.
\end{equation}

A formalism based on the above concepts has been exploited in Ref. \cite{DatZhu} as a tool to work with the matrices describing 
the dynamics of the quarks and neutrinos mixing.  In this context, the matrix $\hat{\eta}$ is recognized as the generatrix of the 
Cabibbo-Kobayashi-Maskawa \cite{ckm} and Pontecorvo-Maki-Nagakawa-Sakata \cite{pmns} matrices\footnote{These matrices 
have complex entries, while all the elements of matrix (\ref{eq:neweta}) are real.}.

In sec. \ref{sec:intro} we have mentioned the addition theorems, which are derived as a consequence of the semi-group property 
of the exponential function. But what about $\e^{\theta\,(h_1 + h_2)}$, with $h_1$ and $h_2$ both satisfying eq. (\ref{eq:genun})? 
Under this hypothesis, taking into account that the ``numbers" $h_{1,2}$ in general do not commute, one has
\begin{equation}
(h_1 + h_2)^2 = a + b_1\,h_1 + b_2\,h_2 + h_1\,h_2 + h_2\,h_1\,,
\end{equation}
that, unless
\begin{equation}
h_1\,h_2 + h_2\,h_1 = b_1\,h_1 + b_2\,h_1\,,
\end{equation}
cannot be put in the form (\ref{eq:genun}). In the 2-dimensional case the problem can be overcome by noting that we can always write
\begin{equation}
\hat{h}_1 + \hat{h}_2 = \hat{1} + \sum_{i = 1}^3 c_i\,\hat{\sigma}_i
\end{equation}
where $\hat{\sigma}_i$ are the Pauli matrices, and therefore we get
\begin{equation}
(\hat{h}_1 + \hat{h}_2)^2 = a + 2\,(\hat{h}_1 + \hat{h}_2) \qquad\qquad a = - 1 + \sum_{i = 1}^3 c_i^2\,,
\end{equation}
i.e., the same structure of eq. (\ref{eq:genun}). This implies that, analogously to eq. (\ref{eq:eulgen}), we can write
\begin{equation}
\e^{\,\theta\,(\hat{h}_1 + \hat{h}_2)} = C (\theta) + (\hat{h}_1 + \hat{h}_2)\,S (\theta)\,.
\end{equation}
For higher order representation we can use different realizations of the Clifford numbers, e.g. the Dirac matrices in the 4-dimensional 
case.

The link of the above procedure to the Hermite-heat polynomials and integral transforms of the Airy type will be discussed in a forthcoming 
investigation, where we will also analyze the properties of ``complex" functions for which the following decomposition holds
\begin{equation}
\label{eq:fcomplex}
f (x + h\,y) = u (x, y) + h\,v (x, y)
\end{equation}
where the functions $u$ and $v$ satisfy the Cauchy-Riemann condition 
\begin{equation}
\partial_y\,\mathbf{w} (x, y) = \hat{h}\,\partial_x\,\mathbf{w} (x, y) \qquad \qquad (\mathbf{w} = (u, v)^T)
\end{equation}
that, taking into account eq. (\ref{eq:derCS}), reflects in the following partial differential equation to be satisfied by the functions $u$ and 
$v$
\begin{equation}
\left(\partial_y^2 - a\,\partial_x^2 - b\,\partial_{xy}^2\right)\,\mathbf{w} = 0\,.
\end{equation}
Therefore, the functions $u$ and $v$ can be considered a generalization of the ordinary harmonic functions.

Furthermore, if the function $f$ admits a Fourier transform, we can write 
\begin{eqnarray}
f (x + h\,y) \speq \frac1{\sqrt{2\,\pi}}\,\int_{- \infty}^\infty \mathrm{d}k\,\e^{\,\imath\,k\,(x + h\,y)}\,f(k) \nn \\
                 \speq \frac1{\sqrt{2\,\pi}}\,\int_{- \infty}^\infty \mathrm{d}k\,\e^{\,\imath\,k\,x}\,f (k)\,
                 \left[C (\imath\,k\,y) + h\,S (\imath\,k\,y)\right]
\end{eqnarray}
thus getting for the functions $u$ and $v$ the following integral representation
\begin{equation}
\left(
\begin{array}{c}
      u (x, y)    \\
      v (x, y)   
\end{array}
\right) = \frac1{\sqrt{2\,\pi}}\,\int_{- \infty}^\infty \mathrm{d}k\,\e^{\,\imath\,k\,x}\,f(k)\,\left(
\begin{array}{c}
      C (\imath\,k\,y)  \\
      S (\imath\,k\,y)
\end{array}
\right)\,.
\end{equation}
The previous relations can be viewed as an extension of the concept of the ordinary Fourier transform.  A further element which can be brought to 
the discussion is the possibility of extending the realm of special functions by introducing new families of Bessel functions. We can indeed use the 
following extension of the Jacobi-Anger generating function
\begin{equation}
\label{eq:jacang}
\e^{\,x\,S (\theta)} = \sum_{m, n = - \infty}^\infty \e^{\,(\alpha\,m + \beta\,n)\,\theta}\,B_{m,n} (x)
\end{equation}
to introduce the two-index Bessel functions $B_{m,n} (x)$, with $S (\theta)$ given by eq. (\ref{eq:cossin}). The recurrence relations satisfied by 
$B_{m,n} (x)$  can easily be obtained from their definition itself. By keeping the derivative of both sides of eq. (\ref{eq:jacang}) with respect to $x$, 
and equating the like-power coefficients in $\e^{\,\alpha\,m\,\theta}$ and $\e^{\,\beta\,n\,\theta}$, we obtain
\begin{equation}
(\alpha - \beta)\,\frac{\mathrm{d}}{\mathrm{d} x}\,B_{m,n} (x) = B_{m - 1,n} (x) - B_{m,n - 1} (x)\,,
\end{equation}
while by deriving with respect to $\theta$ we find
\begin{equation}
(\alpha\,m + \beta\,n)\,B_{m,n} (x) = \frac{x}{\alpha - \beta}\,\left[\alpha\,B_{m - 1,n} (x) - \beta\,B_{m,n - 1} (x)\right]\,.
\end{equation}
It is evident that an analogous function can be associated with the cos-like partner of the $S (\theta)$ function, and that the same procedure 
can be exploited for the higher order trigonometric functions. This is a fairly promising way to improve the theory of Bessel like functions.
These type of problems, along with those relevant to properties of the complex functions (\ref{eq:fcomplex}) will be considered in a forthcoming 
investigation, where we will study more carefully the problem associated with the non-commutativity and the exponential disentanglement.

\end{document}